# FRACTIONAL MARTINGALES AND CHARACTERIZATION OF THE FRACTIONAL BROWNIAN MOTION

By Yaozhong Hu[1], David Nualart[2] and Jian Song

*University of Kansas*

In this paper we introduce the notion of fractional martingale as the fractional derivative of order $\alpha$ of a continuous local martingale, where $\alpha \in (-\frac{1}{2}, \frac{1}{2})$, and we show that it has a nonzero finite variation of order $\frac{2}{1+2\alpha}$, under some integrability assumptions on the quadratic variation of the local martingale. As an application we establish an extension of Lévy's characterization theorem for the fractional Brownian motion.

**1. Introduction.** The fractional Brownian motion (fBm) with Hurst parameter $H \in (0,1)$ is a zero mean Gaussian process with covariance

$$(1.1) \qquad E(B_t^H B_s^H) = \tfrac{1}{2}(t^{2H} + s^{2H} - |t-s|^{2H}).$$

This process is a Brownian motion when $H = \frac{1}{2}$. From the relation $E(|B_t^H - B_s^H|^2) = |t-s|^{2H}$, it follows that $B^H$ has Hölder continuous trajectories of order $H - \varepsilon$, for any $\varepsilon > 0$. On the other hand, the self-similarity of the fBm and the ergodic theorem imply that the fBm has $\frac{1}{H}$-variation on any time interval $[0,t]$ which equals to $c_H t$, where $c_H = E(|B_1^H|^{1/H})$ (see [10]). We refer to the monograph [4] and the review paper [9] for detailed accounts on the properties of the fBm.

In the case of Brownian motion, the famous Lévy's characterization theorem states that a continuous stochastic process $(B_t, t \geq 0)$ adapted to a right-continuous filtration $(\mathcal{F}_t, t \geq 0)$ is an $\mathcal{F}_t$-Brownian motion if and only if $B$ is a local martingale and $\langle B \rangle_t = t$. A natural problem is the extension of Lévy's characterization theorem to the fractional Brownian motion.

Received November 2007; revised February 2009.
[1]Supported by the NSF Grant DMS-05-04783.
[2]Supported by the NSF Grant DMS-06-04207.

The purpose of this paper is to introduce and study the notion of a fractional martingale, and apply it to the above problem. Fix $\alpha \in (-\frac{1}{2}, \frac{1}{2})$. If $M = (M_t, t \geq 0)$ is a continuous local martingale, we denote by $M^{(\alpha)} = (M_t^{(\alpha)}, t \geq 0)$ the stochastic process defined by

$$M_t^{(\alpha)} = \int_0^t (t-s)^\alpha \, dM_s, \tag{1.2}$$

provided this stochastic integral exists for all $t \geq 0$. The process $M^{(\alpha)}$ is called the Riemann–Liouville process of $M$. Notice that $M^{(\alpha)}$ is no longer a martingale and we will say that it is a fractional martingale.

If $\alpha \in (0, \frac{1}{2})$, then the stochastic integral in (1.2) always exists, and $M_t^{(\alpha)} = \Gamma(1+\alpha) I_{0+}^\alpha (M)_t$, where $I_{0+}^\alpha$ is the left-sided fractional integral of order $\alpha$. If $\alpha \in (-\frac{1}{2}, 0)$ and $M$ has $\alpha'$-Hölder continuous trajectories on any finite interval for some $\alpha' > -\alpha$, then $M_t^{(\alpha)}$ exists and $M_t^{(\alpha)} = \Gamma(1+\alpha) D_{0+}^{-\alpha}(M)_t$, where $D_{0+}^{-\alpha}$ is the left-sided fractional derivative of order $-\alpha$. We refer to Samko, Kilbas and Marichev [11] for the definition and properties of the fractional operators.

We are interested in the variation properties of fractional martingales. The process $M^{(\alpha)}$ has Hölder continuous trajectories of order $\gamma$ on any finite interval, for any $\gamma < \frac{1}{2} + \alpha$, provided $M$ has Hölder continuous trajectories of order $\frac{1}{2} - \epsilon$ on any finite interval, for any $\epsilon > 0$. Then, it is natural to expect that $M^{(\alpha)}$ has a finite and nonzero variation of order $\beta = (\frac{1}{2} + \alpha)^{-1} = \frac{2}{1+2\alpha}$. We show that (see Theorem 2.6) if $d\langle M \rangle_t = \xi_t^2 \, dt$, then $M^{(\alpha)}$ has a finite $\beta$-variation $c_\alpha \int_0^t |\xi_s|^\beta \, ds$ under some integrability conditions on $\xi$, where $c_\alpha$ is a constant depending only on $\alpha$. The proof of this result is based on the variation properties of the fractional Brownian motion.

The fractional Brownian motion $B^H$ is not a martingale unless $H = \frac{1}{2}$. But the process

$$M_t = \int_0^t s^{1/2-H}(t-s)^{1/2-H} \, dB_s^H \tag{1.3}$$

is a martingale with respect to the filtration generated by the fBm, verifying $\langle M \rangle_t = d_H t^{2H}$ for some constant $d_H$ (see Norros, Valkeila and Virtamo [8]). We show that if $B = (B_t, t \geq 0)$ is a continuous square integrable centered process with $B_0 = 0$, then $B$ is a fractional Brownian motion with Hurst parameter $H$ if and only if the process $B$ has the following properties:

(i) The sample paths of the process $B$ are Hölder continuous of order $\gamma$ for any $\gamma \in (0, H)$.
(ii) The process $M$ defined in (1.3), where $B^H$ is replaced by $B$, is a martingale with respect to the filtration generated by $B$. If $H > \frac{1}{2}$, we also assume that the quadratic variation of $M$ is absolutely continuous with respect to the Lebesgue measure.

(iii) For any $t > 0$, the process $B$ has $\frac{1}{H}$-variation (in the sense of Definition 2.3) which equals to $c_H t$ on the interval $[0, t]$.

In order to prove that the conditions (i), (ii) and (iii) imply that $B$ is a fractional Brownian motion, it suffices to show that the martingale $M$ satisfies $\langle M \rangle_t = d_H t^{2H}$ for some constant $d_H$, and this will be a consequence of the condition (iii) and the general result on the $\beta$-variation of a fractional martingale.

In a recent work [7], Mishura and Valkeila have proved another extension of the Lévy characterization theorem, where condition (iii) is replaced by an assumption on the renormalized quadratic variation, and no restriction on the quadratic variation of $M$ is required.

THEOREM 1.1 (Mishura and Valkeila). *Assume that $B$ is a continuous square integrable centered process with $B_0 = 0$. Then the following are equivalent:*

(a) *The process $B$ is a fractional Brownian motion with Hurst parameter $H \in (0, 1)$.*
(b) *The process $B$ satisfies the following properties:*
 (i) *The process $B$ has Hölder continuous sample paths of order $\gamma$ for any $\gamma \in (0, H)$ in any finite interval.*
 (ii) *The process $M$ defined in (1.3), where $B^H$ is replaced by $B$, is a martingale with respect to the filtration generated by $B$.*
 (iii) *For any $t > 0$,*
$$\lim_{n \to \infty} n^{2H-1} \sum_{k=1}^{n} (B_{tk/n} - B_{t(k-1)/n})^2 = t^{2H},$$
*in $L^1$.*

The proof of this theorem uses different kind of techniques, and is based on the stochastic calculus with respect to the fractional Brownian motion.

The paper is organized as follows. Section 2 is devoted to study the $\beta$-variation of fractional martingales, and Section 3 contains the proof of the Lévy characterization theorem for the fBm. Some technical lemmas are included in the Appendix.

**2. $\beta$-variation of fractional martingales.** Let $(\Omega, \mathcal{F}, P)$ be a complete probability space equipped with a right-continuous filtration $(\mathcal{F}_t, t \geq 0)$ such that $\mathcal{F}_0$ contains the $P$-null sets. Fix a parameter $\alpha \in (-\frac{1}{2}, \frac{1}{2})$. We introduce the following notion.

DEFINITION 2.1. A continuous $\mathcal{F}_t$-adapted process $(M_t^{(\alpha)}, t \geq 0)$ is called a fractional martingale of order $\alpha$ if there is a continuous local martingale $(M_t, t \geq 0)$ such that, for all $t \geq 0$,

$$\int_0^t (t-s)^{2\alpha} \, d\langle M \rangle_s < \infty, \tag{2.1}$$

almost surely, and

$$M_t^{(\alpha)} = \int_0^t (t-s)^\alpha \, dM_s. \tag{2.2}$$

Notice that by Fubini's theorem condition (2.1) holds true for almost all $t \geq 0$.

If $\alpha \in (0, \frac{1}{2})$, then (2.1) is always fulfilled. Moreover, an integration by parts implies that the integral appearing in (2.2) exists as a Riemann–Stieltjes integral and $M_t^{(\alpha)} = \Gamma(\alpha+1) I_{0+}^\alpha (M)_t$, where $I_{0+}^\alpha$ is the left-sided fractional integral of order $\alpha$.

For any $\alpha \in (-\frac{1}{2}, 0)$ we introduce the following hypothesis:

(H). The trajectories of $M$ are $\alpha'$-Hölder continuous on finite intervals for some $\alpha' > -\alpha$.

Then we have the following result.

LEMMA 2.2. *Fix $\alpha \in (-\frac{1}{2}, 0)$, and let $M$ be a continuous local martingale satisfying condition (H). Then (2.1) holds, $M_t^{(\alpha)}$ exists as a Riemann–Stieltjes integral and it coincides with $\Gamma(\alpha+1) D_{0+}^{-\alpha}(M)_t$, where $D_{0+}^{-\alpha}$ is the left-sided fractional derivative of order $-\alpha$.*

PROOF. Set

$$Z_t = |M_t| + \langle M \rangle_t + \sup_{0 \leq s < u \leq t} \frac{|M_s - M_u|}{|s-u|^{\alpha'}}.$$

For any integer $n \geq 1$ we define

$$T_N = \inf\{t \geq 0 : Z_t > N\}.$$

Then, $T_N$ is an nondecreasing sequence of stopping times such that $T_N \uparrow \infty$. For any $s < t$ we can write

$$E(|\langle M \rangle_{t \wedge T_N} - \langle M \rangle_{s \wedge T_N}|^p) \leq C_p E(|M_{t \wedge T_N} - M_{s \wedge T_N}|^{2p}) \leq C_p N^{2p} |t-s|^{2p\alpha'}.$$

By Kolmogorov's continuity criterion the sample paths of $\langle M \rangle$ are Hölder continuous of order $\gamma$ for any $\gamma < 2\alpha'$, on any finite interval. This implies (2.1), and it is easy to check that the stochastic integral is a Riemann–Stieltjes integral and coincides with $\Gamma(\alpha+1)D_{0+}^{-\alpha}(M)_t$. □

From fractional calculus, assuming condition (H) if $\alpha < 0$, we have $M_t = \frac{1}{\Gamma(\alpha+1)}I_{0+}^{-\alpha}(M^{(\alpha)})_t$, where $I^{-\alpha} = D^\alpha$ if $\alpha > 0$. Using the definition of the left-sided fractional integral and derivative, we have

$$(2.3) \quad M_t = \begin{cases} \dfrac{1}{\Gamma(1+\alpha)\Gamma(-\alpha)} \displaystyle\int_0^t (t-s)^{-1-\alpha} M_s^{(\alpha)} \, ds, & \text{if } \alpha < 0, \\ \dfrac{1}{\Gamma(1+\alpha)\Gamma(1-\alpha)} \displaystyle\int_0^t (t-s)^{-\alpha} \, dM_s^{(\alpha)}, & \text{if } \alpha > 0. \end{cases}$$

In order to define the $\beta$-variation, let us first introduce some notation. Fix a time interval $[a,b]$, and consider the uniform partition

$$\pi^n = \{a = t_0^n < t_1^n < \cdots < t_n^n = b\},$$

where $t_i^n = a + \frac{i}{n}(b-a)$ for $i = 0, \ldots, n$. Let $\beta \geq 1$ and let $X = (X_t, t \geq 0)$ be a continuous stochastic process.

DEFINITION 2.3. We define the $\beta$-variation of $X$ on the interval $[a,b]$, denoted by $\langle X \rangle_{\beta,[a,b]}$, as the limit in probability of

$$(2.4) \quad S_{\beta,n}^{[a,b]}(X) := \sum_{i=1}^n |\Delta_i^n X|^\beta,$$

if the limit exists, where $\Delta_i^n X = X_{t_i^n} - X_{t_{i-1}^n}$. We say that the $\beta$-variation of $X$ on $[a,b]$ exists in $L^1$ if the above limit exists in $L^1$.

We also denote $\langle X \rangle_{\beta,[0,t]}$ by $\langle X \rangle_{\beta,t}$. For instance, a continuous local martingale has a finite 2-variation, denoted by $\langle M \rangle_t$, and the fractional Brownian motion $B_t^H$ of Hurst parameter $H \in (0,1)$ has $\frac{1}{H}$-variation which is equal to $c_H t$, where $c_H = E(|B_1^H|)^{1/H}$.

A direct consequence of the above definition is that if $\langle X \rangle_{\beta,[a,c]}$ exists, then for any $a < b < c$, both $\langle X \rangle_{\beta,[a,b]}$ and $\langle X \rangle_{\beta,[b,c]}$ exist and

$$(2.5) \quad \langle X \rangle_{\beta,[a,c]} = \langle X \rangle_{\beta,[a,b]} + \langle X \rangle_{\beta,[b,c]}.$$

It is also easy to see that the following triangular inequality holds:

$$(2.6) \quad S_{\beta,n}^{[a,b]}(X+Y)^{1/\beta} \leq S_{\beta,n}^{[a,b]}(X)^{1/\beta} + S_{\beta,n}^{[a,b]}(Y)^{1/\beta}.$$

This inequality implies that if $X$ and $Y$ are two continuous stochastic processes such that $\langle X \rangle_{\beta,[a,b]}$ exists and $\langle Y \rangle_{\beta,[a,b]} = 0$, then

$$(2.7) \quad \langle X+Y \rangle_{\beta,[a,b]} = \langle X \rangle_{\beta,[a,b]}.$$

Let $W = (W_t, t \geq 0)$ be an $\mathcal{F}_t$-Brownian motion. We want to compute the $\beta$-variation of $M^{(\alpha)}$, where $M$ is a martingale of the form $M_t = \int_0^t \xi_s \, dW_s$. We will denote by $C$ a generic constant that may depend on $\alpha$. Consider first the case where the martingale is just a standard Wiener process. We recall that

$$\beta = \frac{2}{1+2\alpha}.$$

LEMMA 2.4. *Let $(W_t, t \geq 0)$ be a Wiener process, and set $X_t = W_t^{(\alpha)} = \int_0^t (t-s)^\alpha \, dW_s$. Then the $\beta$-variation of $X$ exists in $L^1$ and $\langle X \rangle_{\beta,[a,b]} = c_\alpha(b-a)$, where $c_\alpha = c_H \kappa_H^{-1/H}$, $H = \frac{1}{2} + \alpha$, $c_H = E(|B_1^H|^{1/H})$, and*

(2.8) $$\kappa_H = \left( \frac{2H\Gamma(3/2 - H)}{\Gamma(H+1/2)\Gamma(2-2H)} \right)^{1/2}.$$

PROOF. Because of (2.5), it is sufficient to show that $\langle X \rangle_{\beta,t} = c_\alpha t$. We can extend the underlying probability space in such a way that $(W_{-t}, t \geq 0)$ is a Brownian motion independent of $W$. Then, the process $B^H$ defined by

$$B_t^H = \kappa_H \left( \int_0^t (t-s)^\alpha \, dW_s + \int_{-\infty}^0 ((t-s)^\alpha - (-s)^\alpha) \, dW_s \right),$$

is a fractional Brownian motion with Hurst parameter $H$ (see Mandelbrot and Van Ness [6]). Hence,

$$X_t = \kappa_H^{-1} B_t^H - Z_t,$$

where $Z_t = \int_{-\infty}^0 ((t-s)^\alpha - (-s)^\alpha) \, dW_s$. From the $\frac{1}{H}$-variation property of fractional Brownian motion we know that $\langle B^H \rangle_{\beta,t} = c_H t$, in $L^1$, because $\beta = \frac{1}{H}$. Then, by (2.7) it suffices to show that $\lim_{n \to \infty} E(|S_{\beta,n}^{[0,t]}(Z)|) = 0$ for all $t \geq 0$. We have

$$\sum_{i=1}^n E(|Z_{t_i^n} - Z_{t_{i-1}^n}|^\beta)$$

$$= C \sum_{i=1}^n \left( \int_{-\infty}^0 ((t_i^n - s)^\alpha - (t_{i-1}^n - s)^\alpha)^2 \, ds \right)^{\beta/2}$$

$$= C \sum_{i=1}^n \left( \int_0^\infty \left( \left( t_{i-1}^n + \frac{t}{n} + s \right)^\alpha - (t_{i-1}^n + s)^\alpha \right)^2 ds \right)^{\beta/2}$$

$$\leq C \left( \int_0^\infty \left( \left( \frac{t}{n} + s \right)^\alpha - s^\alpha \right)^2 ds \right)^{\beta/2} + \frac{C}{n^\beta} \sum_{i=2}^n \left( \int_0^\infty (t_{i-1}^n + s)^{2\alpha - 2} \, ds \right)^{\beta/2}$$

$$= I_1 + I_2.$$

It is easy to see by the dominated convergence theorem that $I_1 \to 0$ as $n \to \infty$. On the other hand,

$$I_2 \leq Ctn^{-1} \sum_{i=2}^{n} (i-1)^{(2\alpha-1)\beta/2} \leq Ctn^{(2\alpha-1)/(2\alpha+1)} \to 0$$

since $\alpha < 1/2$. This proves the lemma. $\square$

We will make use of the following lemma.

LEMMA 2.5. *Fix $a > 0$. For $t \geq a$ let $X_t = \int_0^a (t-s)^\alpha \, dW_s$, where $W = (W_t, t \geq 0)$ is a Wiener process. Then, for all $t \geq a$,*

$$\lim_{n \to \infty} E(|S_{\beta,n}^{[a,t]}(X)|) = 0. \tag{2.9}$$

PROOF. Take $\beta = 2/(1+2\alpha)$. First we have

$$\sum_{i=1}^{n} E \left| \int_0^a [(t_i^n - s)^\alpha - (t_{i-1}^n - s)^\alpha] \, dW_s \right|^\beta$$

$$\leq C \sum_{i=1}^{n} \left\{ \int_0^a [(t_i^n - s)^\alpha - (t_{i-1}^n - s)^\alpha]^2 \, ds \right\}^{\beta/2},$$

where $t \geq a$ and $\{t_i^n\}$ is a uniform partition on $[a,t]$. Then we apply a similar argument as in the proof of Lemma 2.4. $\square$

The following theorem is the main result of this section.

THEOREM 2.6. *Set $\beta = 2/(1+2\alpha)$. Consider a continuous local martingale of the form $M_t = \int_0^t \xi_s \, dW_s$, where $\xi = (\xi_t, t \geq 0)$ is a progressively measurable process such that, for all $t \geq 0$,*

$$\begin{cases} \int_0^t (E(|\xi_s|^\beta))^{\beta'/\beta} \, ds < \infty & \text{for some } \beta' > \beta, \quad \text{if } \alpha < 0, \\ \int_0^t (E(\xi_s^2))^{\beta/2} \, ds < \infty, & \text{if } \alpha > 0. \end{cases} \tag{2.10}$$

*Then, the $\beta$-variation of $M^{(\alpha)}$ on any interval $[0,t]$ exists in $L^1$, and $\langle M^{(\alpha)} \rangle_{\beta,t} = c_\alpha \int_0^t |\xi_s|^\beta \, ds$, where $c_\alpha = c_H \kappa_H^{-1/H}$, $H = \frac{1}{2} + \alpha$, and $\kappa_H$ is defined in (2.8).*

PROOF. We can represent the martingale $M$ as a stochastic integral $M_t = \int_0^t \xi_s \, dW_s$, where $W = (W_t, t \geq 0)$ is a Brownian motion defined on an extension $(\widetilde{\Omega}, \widetilde{\mathcal{F}}, \widetilde{P})$ of our original probability space $(\Omega, \mathcal{F}, P)$. The space $(\widetilde{\Omega}, \widetilde{\mathcal{F}}, \widetilde{P})$ is the product of $(\Omega, \mathcal{F}, P)$, and another space $(\widehat{\Omega}, \widehat{\mathcal{F}}, \widehat{P})$ supporting a Brownian motion independent of $M$. Clearly, if the conclusion of the theorem holds in the extended space, it also holds in the original space.

Notice that if $\alpha < 0$, by Hölder's inequality condition (2.10) implies that

$$\int_0^t (t-s)^{-2\alpha} E(\xi_s^2) \, ds < \infty,$$

and (2.1) holds.

Suppose first that the process $\xi$ has the form $\xi_t = Y I_{(t_1, t_2]}(t)$, where $0 \le t_1 < t_2$ and $Y$ is a bounded $\mathcal{F}_{t_1}$-measurable random variable. In this case the process $M^{(\alpha)}$, denoted by $X$, is given by

$$X_t = Y I_{[t_1, \infty)}(t) \int_{t_1}^{t \wedge t_2} (t-s)^\alpha \, dW_s.$$

For $t \in [0, t_1]$, we clearly have $\langle X \rangle_{\beta, t} = 0$. For $t \in [t_1, t_2]$,

$$X_t = Y \int_0^t (t-s)^\alpha \, dW_s - Y \int_0^{t_1} (t-s)^\alpha \, dW_s,$$

and by Lemmas 2.4 and 2.5, for any interval $[a, b] \subset [t_1, t_2]$, the $\beta$-variation of $X$ exists in $L^1$, and

$$\langle X \rangle_{\beta, [a,b]} = c_\alpha |Y|^\beta (b-a).$$

Finally, by Lemma 2.5, for any interval $[a, b] \subset [t_2, \infty)$, $\langle X \rangle_{\beta, [a,b]} = 0$, in $L^1$. Hence, we have proved that

$$\langle X \rangle_{\beta, t} = c_\alpha |Y|^\beta (t \wedge t_2 - t_1)_+ = c_\alpha \int_0^t |\xi_s|^\beta \, ds.$$

Let us denote by $\mathcal{S}$ the space of step functions of the form

$$\xi_t = \sum_{i=1}^n Y_i I_{(t_{i-1}, t_i]}(t),$$

where $Y_i$ is $\mathcal{F}_{t_{i-1}}$ measurable and bounded, and $0 = t_0 < \cdots < t_n$. For $\xi \in \mathcal{S}$, we have $X_t = \sum_{i=1}^n X_t^i$, where $X_t^i = \int_0^t \xi_t^i (t-s)^\alpha \, dW_s$ and $\xi_t^i = Y_i I_{(t_{i-1}, t_i]}(t)$. From (2.5) we have

$$\langle X \rangle_{\beta, t} = \sum_{i=1}^n \langle X \rangle_{\beta, [t_{i-1}, t_i] \cap [0, t]}.$$

From the first part of the proof we see that

$$\langle X^j \rangle_{\beta, [t_{i-1}, t_i] \cap [0, t]} = \begin{cases} c_\alpha |Y_i|^\beta (t_i \wedge t - t_{i-1})_+, & \text{if } j = i, \\ 0, & \text{if } j \ne i, \end{cases}$$

and applying the triangular inequality (2.6), we see then that

$$\langle X \rangle_{\beta,[t_{i-1},t_i]\cap[0,t]} = \langle X^i \rangle_{\beta,[t_{i-1},t_i]\cap[0,t]}.$$

Hence,

$$(2.11) \qquad \langle X \rangle_{\beta,[0,t]} = c_\alpha \sum_{i=1}^n |Y_i|^\beta (t_i \wedge t - t_{i-1})_+ = c_\alpha \int_0^t |\xi_s|^\beta \, ds,$$

and this proves the result for step functions.

To complete the proof, we use a density argument. Fix a time interval $[0,T]$. We can find a sequence of step functions $(\xi^k, k \geq 1)$ in $\mathcal{S}$ such that if $\alpha > 0$, then

$$\lim_{k\to\infty} \int_0^T (E(|\xi_s - \xi_s^k|^2))^{\beta/2} \, ds = 0,$$

and if $\alpha < 0$, then

$$\lim_{k\to\infty} \int_0^T (E(|\xi_s - \xi_s^k|^\beta))^{\beta'/\beta} \, ds = 0.$$

Define $X_t^k = \int_0^t (t-s)^\alpha \xi_s^k \, dB_s$ for $t \in [0,T]$. From the triangular inequality (2.6) and the Burkholder–Davis–Gundy inequality (see, for instance, [5]), we have, for all $t \in [0,T]$,

$$E(|S_{\beta,n}^{[0,t]}(X)^{1/\beta} - S_{\beta,n}^{[0,t]}(X^k)^{1/\beta}|)$$

$$\leq E((S_{\beta,n}^{[0,t]}(X - X^k))^{1/\beta})$$

$$\leq C \left( E \left( \sum_{i=1}^n \left| \int_0^{t_i^n} ((t_i^n - s)^\alpha - (t_{i-1}^n - s)_+^\alpha)(\xi_s - \xi_s^k) \, dW_s \right|^\beta \right) \right)^{1/\beta}$$

(2.12)

$$\leq C \left( E \left( \sum_{i=1}^n \left| \int_0^{t_i^n} ((t_i^n - s)^\alpha - (t_{i-1}^n - s)_+^\alpha)^2 (\xi_s - \xi_s^k)^2 \, ds \right|^{\beta/2} \right) \right)^{1/\beta}.$$

Now we will consider two cases depending on the sign of $\alpha$.

(i) If $\alpha > 0$, namely, $\beta < 2$, then by the concavity of $x^{\beta/2}$ and Lemma A.1, we have

$$E(|S_{\beta,n}^{[0,t]}(X)^{1/\beta} - S_{\beta,n}^{[0,t]}(X^k)^{1/\beta}|)$$

$$(2.13) \qquad \leq C \left( \sum_{i=1}^n \left| \int_0^{t_i^n} ((t_i^n - s)^\alpha - (t_{i-1}^n - s)_+^\alpha)^2 E(|\xi_s - \xi_s^k|^2) \, ds \right|^{\beta/2} \right)^{1/\beta}$$

$$\leq C \left( \int_0^t (E(|\xi_s - \xi_s^k|^2))^{\beta/2} \, ds \right)^{1/\beta}.$$

Then

$$E\left( \left| S_{\beta,n}^{[0,t]}(X)^{1/\beta} - \left( c_\alpha \int_0^t |\xi_s|^\beta \, ds \right)^{1/\beta} \right| \right)$$

$$\leq E(|S_{\beta,n}^{[0,t]}(X)^{1/\beta} - S_{\beta,n}^{[0,t]}(X^k)^{1/\beta}|)$$

$$+ E\left( \left| S_{\beta,n}^{[0,t]}(X^k)^{1/\beta} - \left( c_\alpha \int_0^t |\xi_s^k|^\beta \, ds \right)^{1/\beta} \right| \right)$$

$$+ c_\alpha^{1/\beta} E\left( \left| \left( \int_0^t |\xi_s^k|^\beta \, ds \right)^{1/\beta} - \left( \int_0^t |\xi_s|^\beta \, ds \right)^{1/\beta} \right| \right).$$

From (2.13) and (2.11) we obtain

$$\limsup_{n \to \infty} E\left( \left| S_{\beta,n}^{[0,t]}(X)^{1/\beta} - \left( c_\alpha \int_0^t |\xi_s|^\beta \, ds \right)^{1/\beta} \right| \right)$$

$$\leq C \left( \int_0^t (E|\xi_s - \xi_s^k|^2)^{\beta/2} \, ds \right)^{1/\beta}$$

$$+ c_\alpha^{1/\beta} E\left( \left| \left( \int_0^t |\xi_s^k|^\beta \, ds \right)^{1/\beta} - \left( \int_0^t |\xi_s|^\beta \, ds \right)^{1/\beta} \right| \right),$$

and letting $k$ tend to zero, we prove the desired result.

(ii) If $\alpha < 0$, namely, $\beta > 2$, then applying the Minkovski inequality in (2.12) and using Lemma A.2, we have

$$E(|S_{\beta,n}^{[0,t]}(X)^{1/\beta} - S_{\beta,n}^{[0,t]}(X^k)^{1/\beta}|)$$

$$\leq C \left( \sum_{i=1}^n \left| \int_0^{t_i^n} ((t_i^n - s)^\alpha - (t_{i-1}^n - s)_+^\alpha)^2 (E(|\xi_s - \xi_s^k|^\beta))^{2/\beta} \, ds \right|^{\beta/2} \right)^{1/\beta}$$

$$\leq C \left( \int_0^t (E|\xi_s - \xi_s^k|^\beta)^{\beta'/\beta} \, ds \right)^{1/\beta'}.$$

Now in the same way as for the case $\alpha > 0$, we can show

$$\lim_{n \to \infty} E\left( \left| S_{\beta,n}^{[0,t]}(X)^{1/\beta} - \left( c_\alpha \int_0^t |\xi_s|^\beta \, ds \right)^{1/\beta} \right| \right) = 0.$$

This proves the theorem. $\square$

REMARK 2.7. If $\alpha > 0$ and $\int_0^t E(\xi_s^2) \, ds < \infty$, then $\int_0^t (E(\xi_s^2))^{\beta/2} \, ds < \infty$, and the $\beta$-variation of the fractional martingale $M^{(\alpha)}$ exists in $L^1$, and $\langle M^{(\alpha)} \rangle_{\beta,t} = c_\alpha \int_0^t |\xi_s|^\beta \, ds$. Using a localization argument, we can prove that this result remains true with the convergence in probability, for any continuous local martingale such that $\langle M \rangle_t = \int_0^t \xi_s^2 \, ds$ for all $t \geq 0$. On the other hand, if $\alpha < 0$ and $\int_0^t E(|\xi_s|^{\beta'}) \, ds < \infty$ for all $t \geq 0$, and for some $\beta' > \beta$, then the $\beta$-variation of the fractional martingale $M^{(\alpha)}$ exists in $L^1$ and $\langle M^{(\alpha)} \rangle_{\beta,t} = c_\alpha \int_0^t |\xi_s|^\beta \, ds$. As a consequence, again by a localization argument, the result remains true with the convergence in probability, for any continuous local martingale such that $\langle M \rangle_t = \int_0^t \xi_s^2 \, ds$, assuming that $\int_0^t |\xi_s|^{\beta'} \, ds < \infty$ almost surely, for all $t \geq 0$, and for some $\beta' > \beta$.

COROLLARY 2.8. *Consider a continuous local martingale $M = (M_t, t \geq 0)$ with $M_0 = 0$ and $\langle M \rangle_t = \int_0^t \xi_s^2 \, ds$, where $\xi = (\xi_t, t \geq 0)$ is a progressively measurable process. Suppose that $M$ satisfies (2.1) for some $\alpha \in (-\frac{1}{2}, \frac{1}{2})$. Then there exists $C > 0$, such that*

$$\liminf_{n \to \infty} E(S_{\beta,n}^{[a,b]}(M^{(\alpha)})) \geq C \int_a^b E(|\xi_s|^\beta) \, ds.$$

PROOF. For each integer $N \geq 1$ let $\psi_N(x) = x$ if $|x| \leq N$ and $\psi_N(x) = \frac{N}{x}$ if $|x| > N$. Denote $M_t^{(\alpha),N} = \int_0^t (t-s)^\alpha \psi_N(\xi_s) \, dM_s$. An application of Burkholder's inequality yields

$$E(S_{\beta,n}^{[a,b]}(M^{(\alpha)})) = E\left(\sum_{i=1}^n \left| \int_0^{t_i^n} ((t_i^n - s)^\alpha - (t_{i-1}^n - s)_+^\alpha) \, dM_s \right|^\beta \right)$$

$$\geq CE\left(\sum_{i=1}^n \left| \int_0^{t_i^n} ((t_i^n - s)^\alpha - (t_{i-1}^n - s)_+^\alpha)^2 |\xi_s|^2 \, ds \right|^{\beta/2} \right)$$

$$\geq CE\left(\sum_{i=1}^n \left| \int_0^{t_i^n} ((t_i^n - s)^\alpha - (t_{i-1}^n - s)_+^\alpha)^2 (|\xi_s| \wedge N)^2 \, ds \right|^{\beta/2} \right)$$

$$\geq CE(S_{\beta,n}^{[a,b]}(M^{(\alpha),N})).$$

By Theorem 2.6, $S_{\beta,n}^{[a,b]}(M^{(\alpha),N})$ converges to $\int_a^b (|\xi_s| \wedge N)^\beta \, ds$ in $L^1$ as $n$ tends to infinity. So, $\lim_{n \to \infty} E(S_{\beta,n}^{[a,b]}(M^{(\alpha),N})) = \int_a^b E((|\xi_s| \wedge N)^\beta) \, ds$ and, consequently, $\liminf_{n \to \infty} E(S_{\beta,n}^{[a,b]}(M^{(\alpha)})) \geq C \int_a^b E|\xi_s|^\beta \, ds$. □

So far we have considered continuous local martingales such that $\langle M \rangle_t$ is absolutely continuous with respect to the Lebesgue measure. The next result says that in the case $\alpha < 0$ if the quadratic variation of the martingale is not absolutely continuous with respect to the Lebesgue measure with positive probability, then the $\beta$-variation is infinite.

PROPOSITION 2.9. *Fix $-\frac{1}{2} < \alpha < 0$. Suppose that $M = (M_t, t \geq 0)$ is a continuous local martingale, satisfying (2.1). Consider the Lebesgue decomposition of its quadratic variation given by $\langle M \rangle_t = \mu_t + \nu_t$, where $\mu_t$ and $\nu_t$ are continuous nondecreasing adapted processes such that $d\mu_t$ is absolutely continuous with respect to the Lebesgue measure, and $d\nu_t$ is singular. If $P(d\nu_t \neq 0) > 0$, then we have $\lim_{n\to\infty} E(S_{\beta,n}^{[0,t]}(M^{(\alpha)})) = \infty$, for all $t \geq 0$.*

PROOF. By Burkholder's inequality, we have

$$E\left(\sum_{i=1}^{n} |M_{t_i^n}^{(\alpha)} - M_{t_{i-1}^n}^{(\alpha)}|^\beta\right)$$
$$\geq C \sum_{i=1}^{n} E\left(\int_0^{t_i^n} ((t_i^n - s)^\alpha - (t_{i-1}^n - s)_+^\alpha)^2 \, d\langle M\rangle_s\right)^{\beta/2}$$
$$\geq C \sum_{i=1}^{n} E\left(\int_0^{t_i^n} ((t_i^n - s)^\alpha - (t_{i-1}^n - s)_+^\alpha)^2 \, d\mu_s\right)^{\beta/2}$$
$$+ C \sum_{i=1}^{n} E\left(\int_0^{t_i^n} ((t_i^n - s)^\alpha - (t_{i-1}^n - s)_+^\alpha)^2 \, d\nu_s\right)^{\beta/2}.$$

Then the result follows from the above inequality and Lemma A.3, proved in the Appendix. □

On the other hand, the next result says that in the case $\alpha \in (0, \frac{1}{4})$, the $\beta$-variation is zero if the quadratic variation of the martingale is singular.

PROPOSITION 2.10. *Suppose that $M = (M_t, t \geq 0)$ is a continuous local martingale, such that almost surely the measure $d\langle M\rangle_t$ is singular with respect to the Lebesgue measure. Then, if $\alpha \in (0, \frac{1}{4})$, we have $\lim_{n\to\infty} E(S_{\beta,n}^{[0,t]}(M^{(\alpha)})) = 0$, for all $t \geq 0$.*

PROOF. The result is an immediate consequence of Lemma A.3, proved in the Appendix. □

**3. Characterization of fractional Brownian motion.** Suppose that $B^H$ is a fractional Brownian motion with Hurst parameter $H \in (0,1)$. The process $B^H$ admits the following representation (see [4]):

(3.1) $$B_t^H = \int_0^t Z_H(t, s) \, dW_s,$$

where

$$Z_H(t, s) = \kappa_H \left[ \left(\frac{t}{s}\right)^{H-1/2} (t-s)^{H-1/2} \right.$$

(3.2)
$$-\left(H-\frac{1}{2}\right)s^{1/2-H}\int_s^t u^{H-3/2}(u-s)^{H-1/2}\,du\bigg],$$

with $\kappa_H$ defined in (2.8).

The next theorem is the main result of this paper and provides an extension of Lévy characterization to the fractional Brownian motion.

THEOREM 3.1. *Fix $H \in (0,1)$, $H \neq \frac{1}{2}$. Suppose that $B = (B_t, t \geq 0)$ is a zero mean continuous stochastic process. The following two conditions are equivalent:*

(1) *$B$ is a fractional Brownian motion with Hurst parameter $H$.*
(2) *The process $B$ satisfies the following conditions:*
   (i) *The trajectories of $B$ are Hölder continuous of order $H - \epsilon$ for any $H - \epsilon \in (0, H)$.*
   (ii) *Let*

(3.3)
$$M_t = \int_0^t s^{1/2-H}(t-s)^{1/2-H}\,dB_s.$$

   *Then $M$ is a local martingale. Furthermore, if $H > \frac{1}{2}$, the quadratic variation of the martingale $M$ is absolutely continuous with respect to the Lebesgue measure almost surely.*
   (iii) *For any $t > 0$, the $\frac{1}{H}$-variation of $B$ in the interval $[0, t]$ exists in $L^1$, and $\langle B \rangle_{1/H,t} = c_H t$, where $c_H = E(|\xi|^{1/H})$ and $\xi$ is a standard normal random variable.*

REMARK 3.2. Notice that condition (i) is always true if $H < \frac{1}{2}$, and the Riemann–Stieltjes integral in (3.3) exists by Proposition A.6.

PROOF OF THEOREM 3.1. From the properties of the fractional Brownian motion we know that (1) implies (2). Suppose that (2) holds. Fix $H - \epsilon \in (0, H)$, and $T > 0$. We are going to show that $B$ is a fractional Brownian motion with Hurst parameter $H$ in the time interval $[0, T]$. Denote by $\|B\|_{H-\epsilon}$ the Hölder norm of order $H - \epsilon$ on $[0, T]$ [see (A.2)]. The proof is divided into several steps.

*Step* 1. From (3.3), we can solve the integral equation to express $B$ as a functional of $M$. This can be done as in the proof of Theorem 5.2 of [8]. In this way we obtain

$$B_t = d_H[t^{H-1/2}R_t - (H - \tfrac{1}{2})Y_t],$$

where $d_H = B(\tfrac{3}{2} - H, H + \tfrac{1}{2})^{-1}$,

$$R_t = \int_0^t (t-s)^{H-1/2}\,dM_s,$$

and
$$Y_t = \int_0^t \left( \int_s^t u^{H-3/2}(u-s)^{H-1/2} \, du \right) dM_s.$$

Comparing with the representation formula (3.1) for the fractional Brownian motion, it suffices to prove that

(3.4) $$d\langle M \rangle_s = (\kappa_H d_H^{-1} s^{1/2-H})^2 \, ds,$$

because this implies that $M$ is a Gaussian martingale, and $B$ has the covariance of the fractional Brownian motion with Hurst parameter $H$. In order to show (3.4), we are going to compute the $\frac{1}{H}$-variation of $R$, from the decomposition

(3.5) $$R_t = d_H^{-1} t^{1/2-H} B_t + (H - \tfrac{1}{2}) t^{1/2-H} Y_t.$$

*Step* 2. Fix $0 < \epsilon < H \wedge \tfrac{1}{2} \wedge (1-H)$ and suppose that $E(\|B\|_{H-\varepsilon}^{1/H}) < \infty$. We will first show that the $\frac{1}{H}$-variation of the process $Z_t = t^{1/2-H} B_t$ exists in $L^1$ in any interval $[0,t] \subset [0,T]$, and

(3.6) $$\langle Z \rangle_{1/H,t} = 2H c_H t^{1/(2H)}.$$

An application of the triangular inequality yields

(3.7) $$S_{1/H,n}^{[0,t]}(Z) \leq \left| \left( \sum_{i=1}^n (t_i^n)^{1/(2H)-1} |B_{t_i^n} - B_{t_{i-1}^n}|^{1/H} \right)^H \right.$$
$$\left. + \left( \sum_{i=1}^n |(t_i^n)^{1/2-H} - (t_{i-1}^n)^{1/2-H}|^{1/H} |B_{t_{i-1}^n}|^{1/H} \right)^H \right|^{1/H},$$

and

(3.8) $$S_{1/H,n}^{[0,t]}(Z) \geq \left| \left( \sum_{i=1}^n (t_i^n)^{1/(2H)-1} |B_{t_i^n} - B_{t_{i-1}^n}|^{1/H} \right)^H \right.$$
$$\left. - \left( \sum_{i=1}^n |(t_i^n)^{1/2-H} - (t_{i-1}^n)^{1/2-H}|^{1/H} |B_{t_{i-1}^n}|^{1/H} \right)^H \right|^{1/H}.$$

We have
$$\sum_{i=1}^n |(t_i^n)^{1/2-H} - (t_{i-1}^n)^{1/2-H}|^{1/H} |B_{t_{i-1}^n}|^{1/H}$$

(3.9) $$\leq C \|B\|_{H-\varepsilon}^{1/H} \left( \frac{t}{n} \right)^{1/(2H)-\varepsilon/H} \sum_{i=2}^n (i-1)^{-1/(2H)-\varepsilon/H}$$
$$\leq C \|B\|_{H-\varepsilon}^{1/H} t^{1/(2H)-\varepsilon/H} n^{1-1/H},$$

which converges in $L^1$ to 0 as $n$ tends to infinity. From (3.7) to (3.9) we obtain

$$(3.10) \quad \lim_{n\to\infty} S^{[0,t]}_{1/H,n}(Z) = \lim_{n\to\infty} \sum_{i=1}^{n} (t_i^n)^{1/(2H)-1} |B_{t_i^n} - B_{t_{i-1}^n}|^{1/H},$$

in $L^1$, provided that the limit on the right-hand side of (3.10) exists. Denote $I_j^n = (t_{j-1}^n, t_j^n]$ for $j = 1, 2, \ldots, n$. We divide every subinterval $I_j^n$ into $m$ parts, and we get a finer partition $0 = t_0^{nm} < \cdots < t_{nm}^{nm} = t$. Then, we have

$$\left| \sum_{i=1}^{nm} (t_i^{nm})^{1/(2H)-1} |B_{t_i^{nm}} - B_{t_{i-1}^{nm}}|^{1/H} - \sum_{j=1}^{n} c_H (t_j^n)^{1/(2H)-1} (t_j^n - t_{j-1}^n) \right|$$

$$= \left| \sum_{j=1}^{n} \left( \sum_{i=(j-1)m+1}^{jm} ((t_i^{nm})^{1/(2H)-1} - (t_j^n)^{1/(2H)-1}) |B_{t_i^{nm}} - B_{t_{i-1}^{nm}}|^{1/H} \right. \right.$$

$$\left. \left. + (t_j^n)^{1/(2H)-1} \left( \sum_{i=(j-1)m+1}^{jm} |B_{t_i^{nm}} - B_{t_{i-1}^{nm}}|^{1/H} - c_H (t_j^n - t_{j-1}^n) \right) \right) \right|$$

$$\leq \sum_{j=1}^{n} |(t_j^n)^{1/(2H)-1} - (t_{j-1}^n)^{1/(2H)-1}| \sum_{i=(j-1)m+1}^{jm} |B_{t_i^{nm}} - B_{t_{i-1}^{nm}}|^{1/H}$$

$$+ (t_j^n)^{1/(2H)-1} \left| \sum_{i=(j-1)m+1}^{jn} |B_{t_i^{nm}} - B_{t_{i-1}^{nm}}|^{1/H} - c_H (t_j^n - t_{j-1}^n) \right|.$$

Letting $m$ tend to infinity and using assumption (ii), we obtain

$$\lim_{n\to\infty} \sum_{i=1}^{n} (t_i^n)^{1/(2H)-1} |B_{t_i^n} - B_{t_{i-1}^n}|^{1/H} = 2H c_H t^{1/(2H)},$$

in $L^1$, which shows (3.6).

*Step* 3. We claim that the $\frac{1}{H}$-variation of the process $V_t = t^{1/2-H} Y_t$ in $L^1$ is zero. The increment $|Y_t - Y_s|$ can be estimated by Lemma A.7 in the Appendix with $\alpha = \frac{1}{2} - H$, $f$ being a trajectory of the process $B$ and $\beta = H - \varepsilon$. Notice that $\alpha + \beta = \frac{1}{2} - \varepsilon$, and $2\alpha + \beta = 1 - H - \varepsilon$. Hence, for any $s, t \in [0, T]$, we have

$$|Y_t - Y_s| \leq C \|B\|_{H-\varepsilon} (t^\beta - s^\beta).$$

Therefore, as in (3.7), we have

$$E(S^{[0,t]}_{1/H,n}(V)) \leq C \sum_{i=1}^{n} (t_i^n)^{1/(2H)-1} E(|Y_{t_i^n} - Y_{t_{i-1}^n}|^{1/H})$$

$$+ C \sum_{i=1}^{n} ((t_i^n)^{1/2-H} - (t_{i-1}^n)^{1/2-H})^{1/H} E(|Y_{t_{i-1}^n}|^{1/H})$$

$$= A_n + B_n.$$

For the term $A_n$ we have

$$A_n \leq C\|B\|_{H-\varepsilon}^{1/H} \sum_{i=1}^{n} (t_i^n)^{1/(2H)-1} ((t_i^n)^{H-\varepsilon} - (t_{i-1}^n)^{H-\varepsilon})^{1/H}$$

$$= C\|B\|_{H-\varepsilon}^{1/H} \left(\frac{t}{n}\right)^{1/(2H)-\varepsilon/H} \sum_{i=1}^{n} i^{1/(2H)-1} (i-1)^{1-\varepsilon/H-1/H}$$

$$\leq C\|B\|_{H-\varepsilon}^{1/H} \left(\frac{t}{n}\right)^{1/(2H)-\varepsilon/H} n^{-1/(2H)-\varepsilon/H+1}.$$

By Lemma A.7, $\lim_{n\to\infty} E(A_n) = 0$. For the term $B_n$, using that $E(|Y_{t_{i-1}^n}|^{1/H}) \leq CE(\|B\|_{H-\varepsilon}^{1/H})|t_{i-1}^n|^{1-\varepsilon/H}$, we obtain

$$E(B_n) \leq CE(\|B\|_{H-\varepsilon}^{1/H}) \sum_{i=1}^{n} (t_{i-1}^n)^{-1/(2H)-\varepsilon/H} \left(\frac{t}{n}\right)^{1/H}$$

$$\leq CE(\|B\|_{H-\varepsilon}^{1/H}) \left(\frac{1}{n}\right)^{-1+1/H-\varepsilon/H} \to 0.$$

Hence, $\langle Y \rangle_{1/H,t} = 0$, in $L^1$, for all $t \in [0, T]$.

*Step* 4. From (3.5), (3.6), Step 3 and (2.7), we get that the $\frac{1}{H}$-variation of the process $R$ in any interval $[0, t] \subset [0, T]$ exists in $L^1$, and

(3.11) $$\langle R \rangle_{1/H,t} = c_H d_H^{-1/H} 2H t^{1/(2H)}.$$

On the other hand, since $R_t$ is an $H - \frac{1}{2}$ martingale, Theorem 2.6 and Proposition 2.9 imply that if $H < 1/2$, the quadratic variation $d\langle M \rangle_s$ must be absolutely continuous with respect to the Lebesgue measure, almost surely. In the case $H > \frac{1}{2}$ this is true by the assumption (ii). This implies that $\langle M \rangle_t = \int_0^t \xi_s^2 \, ds$, where $\xi = (\xi_t, t \geq 0)$ is a progressively measurable process.

By Corollary 2.8, there is a positive constant $C$ such that, for any $t_1, t_2 \in [0, T]$, $C \int_{t_1}^{t_2} s^{1/(2H)-1} \, ds \geq \int_{t_1}^{t_2} E(|\xi_s|^{1/H}) \, ds$. Then $E(|\xi_s|^{1/H}) \leq Cs^{1/(2H)-1}$. Thus, we can apply Theorem 2.6 to obtain $\langle R \rangle_{1/H,t} = c_H \kappa_H^{-1/H} \int_0^t |\xi_s|^{1/H} \, ds$. Comparing this with (3.11), we obtain

$$|\xi_s| = \kappa_H d_H^{-1} s^{1/2-H}, \qquad 0 \leq s \leq t,$$

and (3.4) holds. This proves that $B$ is a fractional Brownian motion with Hurst parameter $H$ under the condition $E(\|B\|_{H-\varepsilon}^{1/H}) < \infty$.

*Step* 5. If $E(\|B\|_{H-\varepsilon}^{1/H})$ is not necessarily finite, we can use a localization argument. Denote

$$T_K = \inf\{t \geq 0 : \|B\|_{t,H-\varepsilon} \geq K\} \wedge T,$$

and $B_t^K = B_{t \wedge T_K}$. Since $\sum_{i=1}^n |B_{t_i^n}^K - B_{t_{i-1}^n}^K|^{1/H} \leq \sum_{i=1}^n |B_{t_i^n} - B_{t_{i-1}^n}|^{1/H} + (K\frac{t}{n})^{1/H}$, by the dominated convergence theorem, we can also get

$$\lim_n E\left(\left|\sum_{i=1}^n |B_{t_i^n}^K - B_{t_{i-1}^n}^K|^{1/H} - c_H(t \wedge T_K)\right|\right) = 0.$$

By modifying the proof in Steps 1–4 slightly, we get

$$|\xi_s| = \kappa_H d_H^{-1} s^{1/2-H}, \qquad 0 \leq s \leq t \wedge T_K.$$

Clearly, $\lim_{K \to \infty} T_K = T$, and then

$$|\xi_s| = \kappa_H d_H^{-1} s^{1/2-H}, \qquad 0 \leq s \leq T. \qquad \square$$

REMARK 3.3. Notice that in the case $H > \frac{1}{2}$ we have imposed the additional assumption that the martingale (3.3) has an absolutely continuous quadratic variation. This is true, for instance, if the filtration generated by the process $B$ is included in the filtration generated by a Brownian motion. The next proposition shows that this condition is necessary at least in the case $H \in (\frac{1}{2}, \frac{3}{4})$.

PROPOSITION 3.4. *Suppose that $H \in (\frac{1}{2}, \frac{3}{4})$. There exists a process $B$, satisfying conditions* (i) *and* (iii) *of Theorem 3.1, such that the process $M$ defined in (3.3) is a local martingale, and $B$ is not a fractional Brownian motion.*

PROOF. Let $B^H$ be a fractional Brownian motion with Hurst parameter $H \in (\frac{1}{2}, \frac{3}{4})$. Define

$$M_t = \int_0^t s^{1/2-H}(t-s)^{1/2-H} \, dB_s^H.$$

Let $N_t = W_{\phi(t)}$, where $W$ is a Brownian motion independent of $B^H$, and $\phi$ is a strictly increasing, Hölder continuous function of exponent $\gamma$ for any $\gamma < 1$, null at zero, such that the measure $d\phi(t)$ is singular with respect to the Lebesgue measure (for the existence of such function, see Lemma A.8 in the Appendix). Set

$$\widetilde{M}_t = M_t + N_t \quad \text{and} \quad \widetilde{B}_t^H = B_t^H + Y_t,$$

where

$$Y_t = d_H \bigg( t^{H-1/2} \int_0^t (t-s)^{H-1/2} \, dN_s$$
$$- \bigg( H - \frac{1}{2} \bigg) \int_0^t \bigg( \int_s^t u^{H-3/2}(u-s)^{H-1/2} \, du \bigg) dN_s \bigg).$$

The process $\widetilde{B}^H$ clearly satisfies (i) and it is not a fractional Brownian motion. Finally, $\langle \widetilde{B}^H \rangle_{1/H,t} = c_H t$ in $L^1$, because the $\frac{1}{H}$-variation of $\int_0^t (t-s)^{H-1/2} \, dN_s$ is zero by Proposition 2.10, and, by the same arguments as in the proof of Theorem 3.1, we can show that the $\frac{1}{H}$-variation of $Y$ vanishes. □

## APPENDIX

### A.1. Some technical lemmas.

LEMMA A.1. *Let $\alpha \in (0, \frac{1}{2})$. Fix an interval $[0, t]$. For any natural number $m$, we define $t_i^m = \frac{i}{m} t$, $0 \le i \le m$. Let $g$ be a measurable function on $[0, \infty)$ such that, for all $t \ge 0$, $\int_0^t |g(s)| \, ds < \infty$. Then there exists a function $C(t) > 0$ satisfying*

$$\limsup_{m \to \infty} \sum_{i=1}^m \bigg( \int_0^{t_i^m} ((t_i^m - s)^\alpha - (t_{i-1}^m - s)_+^\alpha)^2 |g(s)| \, ds \bigg)^{\beta/2}$$
$$\le C(t) \int_0^t |g(s)|^{\beta/2} \, ds.$$

PROOF. Set

$$A_m = \sum_{i=1}^m \bigg( \int_0^{t_i^m} ((t_i^m - s)^\alpha - (t_{i-1}^m - s)_+^\alpha)^2 |g(s)| \, ds \bigg)^{\beta/2}.$$

We have $A_m \le C(A_{1,m} + A_{2,m} + A_{3,m})$, where

$$A_{1,m} = \sum_{i=3}^m \bigg( \int_0^{t_{i-2}^m} ((t_i^m - s)^\alpha - (t_{i-1}^m - s)^\alpha)^2 |g(s)| \, ds \bigg)^{\beta/2},$$

$$A_{2,m} = \sum_{i=2}^m \bigg( \int_{t_{i-2}^m}^{t_{i-1}^m} ((t_i^m - s)^\alpha - (t_{i-1}^m - s)^\alpha)^2 |g(s)| \, ds \bigg)^{\beta/2}$$

and

$$A_{3,m} = \sum_{i=1}^m \bigg( \int_{t_{i-1}^m}^{t_i^m} (t_i^m - s)^{2\alpha} |g(s)| \, ds \bigg)^{\beta/2}.$$

Let $\phi_m(x) = ((x + \frac{t}{m})^\alpha - x^\alpha)^2$. The $\phi(x)$ is a nonincreasing of $x$ when $x \geq 0$. As a consequence,

$$A_{1,m} = \frac{m}{t} \sum_{i=3}^{m} \int_{t_{i-2}^m}^{t_{i-1}^m} \left( \int_0^{t_{i-2}^m} ((t_i^m - s)^\alpha - (t_{i-1}^m - s)^\alpha)^2 |g(s)| \, ds \right)^{\beta/2} du$$

$$= \frac{m}{t} \sum_{i=3}^{m} \int_{t_{i-2}^m}^{t_{i-1}^m} \left( \int_0^{t_{i-2}^m} \phi_m(t_{i-1}^m - s)|g(s)| \, ds \right)^{\beta/2} du$$

$$\leq \frac{m}{t} \int_0^t \left( \int_0^u \phi_m(u - s)|g(s)| \, ds \right)^{\beta/2} du.$$

Using the Hölder inequality, we obtain

$$\left( \int_0^u \phi_m(u-s)|g(s)|\,ds \right)^{\beta/2} \leq \left( \int_0^u \phi_m(u-s)\,ds \right)^{\beta/2-1}$$
$$\times \int_0^u \phi_m(u-s)|g(s)|^{\beta/2}\,ds$$
$$\leq \left( \int_0^t \phi_m(s)\,ds \right)^{\beta/2-1} \int_0^u \phi_m(u-s)|g(s)|^{\beta/2}\,ds.$$

Integrating in the variable $u$ yields

$$A_{1,m} \leq \frac{m}{t} \left( \int_0^t \phi_m(s)\,ds \right)^{\beta/2-1} \int_0^t \int_0^u \phi_m(u-s)\,ds |g(s)|^{\beta/2}\,ds\,du$$

$$= \frac{m}{t} \left( \int_0^t \phi_m(s)\,ds \right)^{\beta/2-1} \int_0^t \left( \int_s^t \phi_m(u-s)\,du \right) |g(s)|^{\beta/2}\,ds$$

(A.1)

$$\leq \frac{m}{t} \left( \int_0^t \phi_m(s)\,ds \right)^{\beta/2-1} \int_0^t \left( \int_0^t \phi_m(u)\,du \right) |g(s)|^{\beta/2}\,ds$$

$$= \frac{m}{t} \left( \int_0^t \phi_m(s)\,ds \right)^{\beta/2} \int_0^t |g(s)|^{\beta/2}\,ds.$$

Therefore,

$$\lim_{m \to \infty} A_{1,m} = t^{-1} \left( \int_0^\infty ((x+t)^\alpha - x^\alpha)^2 \, dx \right)^{\beta/2} \int_0^t |g(s)|^{\beta/2}\,du.$$

For the term $A_{3,m}$ we can write

$$A_{3,m} \leq \left( \frac{t}{m} \right)^{\alpha\beta} \sum_{i=1}^m \left( \int_{t_{i-1}^m}^{t_i^m} |g(s)|\,ds \right)^{\beta/2}$$

$$= \sum_{i=1}^m \left( \frac{m}{t} \int_{t_{i-1}^m}^{t_i^m} |g(s)|\,ds \right)^{\beta/2} \frac{t}{m}.$$

The functions

$$g_m(s) = \frac{m}{t} \sum_{i=1}^m \left( \int_{t_{i-1}^m}^{t_i^m} |g(s)|\,ds \right) I_{(t_{i-1}^m, t_i^m]}(s)$$

converge almost everywhere to $|g|$, and they are bounded in $L^1([0,t])$. Hence, $|g(s)|^{\beta/2}$ is uniformly integrable on $[0,t]$. Therefore,

$$\limsup_{m\to\infty} A_{3,m} \leq \lim_{m\to\infty} \int_0^t |g_m(s)|^{\beta/2}\,ds = \int_0^t |g(s)|^{\beta/2}\,ds.$$

From the fact that $|x^\alpha - y^\alpha| \leq |x-y|^\alpha$, we see that

$$A_{2,m} \leq \sum_{i=2}^m \left( \int_{t_{i-2}^m}^{t_{i-1}^m} |t_i^m - t_{i-1}^m|^{2\alpha} |g(s)|\,ds \right)^{\beta/2}.$$

Thus, in the same way as for $A_{3,m}$, we have

$$\limsup_{m\to\infty} A_{2,m} \leq 2\int_0^t |g(s)|^{\beta/2}\,ds. \qquad \square$$

LEMMA A.2. *Let $\alpha \in (-\frac{1}{2}, 0)$. Fix an interval $[0,t]$. For any natural number $m$, we define $t_i^m = \frac{i}{m}t$, $0 \leq i \leq m$. Let $g$ be a measurable function on $[0,\infty)$ such that, for all $t \geq 0$, $\int_0^t |g(s)|^{\beta'/2}\,ds < \infty$ for some $\beta' > \beta$. Then there exists a constant $C$ depending on $t$ such that*

$$\sum_{i=1}^m \left( \int_0^{t_i^m} ((t_i^m - s)^\alpha - (t_{i-1}^m - s)_+^\alpha)^2 |g(s)|\,ds \right)^{\beta/2}$$

$$\leq C \left( \int_0^t |g(s)|^{\beta'/2}\,ds \right)^{\beta/\beta'}.$$

PROOF. Consider the decomposition given in the proof of Lemma A.1. For the first term we can write, from inequality (A.1),

$$A_{1,m} \leq C\frac{m}{t} \left( \int_0^t \left( s^{2\alpha} - \left(s + \frac{t}{m}\right)^{2\alpha} \right) ds \right)^{\beta/2} \int_0^t |g(s)|^{\beta/2}\,ds$$

$$\leq C\frac{m}{t} \left( t^{1+2\alpha} + \left(\frac{t}{m}\right)^{1+2\alpha} - \left(t + \frac{t}{m}\right)^{1+2\alpha} \right)^{\beta/2} \int_0^t |g(s)|^{\beta/2}\,ds$$

$$\leq C\frac{m}{t} \left(\frac{t}{m}\right)^{(1+2\alpha)\beta/2} \int_0^t |g(s)|^{\beta/2}\,ds$$

$$\leq C \int_0^t |g(s)|^{\beta/2}\,ds \leq C \left( \int_0^t |g(s)|^{\beta'/2}\,ds \right)^{\beta'/\beta}.$$

Let $2\alpha p > -1$ and $\frac{1}{p} + \frac{1}{q} = 1$. Then $\beta' = 2q > \beta$, and applying Hölder's inequality, we can write

$$A_{3,m} \leq \sum_{i=1}^{m} \left( \int_{t_{i-1}^m}^{t_i^m} (t_i^m - s)^{2\alpha p} \, ds \right)^{\beta/(2p)} \left( \int_{t_{i-1}^m}^{t_i^m} |g(s)|^q \, ds \right)^{\beta/(2q)}$$

$$\leq C \sum_{i=1}^{m} \left( \frac{t}{m} \right)^{((1+2\alpha p)/p)\beta/2} \left( \int_{t_{i-1}^m}^{t_i^m} |g(s)|^q \, ds \right)^{\beta/(2q)}$$

$$\leq C t^{((1+2\alpha p)/p)\beta/2} \left( \int_0^t |g(s)|^q \, ds \right)^{\beta/(2q)}.$$

For the term $A_{2,m}$, with the same notation as above, we can write

$$A_{2,m} \leq C \sum_{i=2}^{m} \left( \int_{t_{i-2}^m}^{t_{i-1}^m} (t_{i-1}^m - s)^{2\alpha} |g(s)| \, ds \right)^{\beta/2}$$

$$\leq C \sum_{i=2}^{m} \left( \frac{t}{m} \right)^{((1+2\alpha p)/p)\beta/2} \left( \int_{t_{i-2}^m}^{t_{i-1}^m} |g(s)|^q \, ds \right)^{\beta/(2q)}$$

$$\leq C t^{((1+2\alpha p)/p)\beta/2} \left( \int_0^t |g(s)|^q \, ds \right)^{\beta/(2q)}. \qquad \square$$

LEMMA A.3. *Suppose that $v$ is a measure on an interval $[0,t]$, which is singular with respect to the Lebesgue measure. We have the following:*

(i) *If $\alpha \in (-\frac{1}{2}, 0)$, then*

$$\lim_{n \to \infty} \sum_{i=1}^{n} \left( \int_0^{t_i^n} ((t_i^n - s)^\alpha - (t_{i-1}^n - s)_+^\alpha)^2 \, d\nu_s \right)^{\beta/2} = \infty.$$

(ii) *If $\alpha \in (0, \frac{1}{4})$, then*

$$\lim_{n \to \infty} \sum_{i=1}^{n} \left( \int_0^{t_i^n} ((t_i^n - s)^\alpha - (t_{i-1}^n - s)_+^\alpha)^2 \, d\nu_s \right)^{\beta/2} = 0.$$

PROOF. Denote $\triangle_i^n := (t_{i-1}^n, t_i^n]$. Set

$$A_n = \sum_{i=1}^{n} \left( \int_0^{t_i^n} ((t_i^n - s)^\alpha - (t_{i-1}^n - s)_+^\alpha)^2 \, d\nu_s \right)^{\beta/2}.$$

(i) If $\alpha \in (-\frac{1}{2}, 0)$, then

$$A_n \geq \sum_{i=1}^{n} \left( \int_{t_{i-1}^n}^{t_i^n} (t_i^n - s)^{2\alpha} \, d\nu_s \right)^{\beta/2}$$

$$\geq C\left(\frac{t}{n}\right)^{\alpha\beta}\sum_{i=1}^{n}(\nu(\triangle_i^n))^{\beta/2} \geq \sum_{i=1}^{n}C\left(\frac{t}{n}\right)\left(\frac{\nu(\triangle_i^n)}{m(\triangle_i^n)}\right)^{\beta/2},$$

where $m$ denotes the Lebesgue measure. Suppose that $\mathcal{F}_n$ is the $\sigma$-field of subsets of the interval $[0,t]$ generated by the partition $\{\triangle_i^n, i=1,\ldots,n\}$. Denote by $\nu_n$ and $m_n$ the restrictions of the measures $\nu$ and $m$ to the $\sigma$-field $\mathcal{F}_n$. Set

$$X_n = \sum_{i=1}^{n}\frac{\nu(\triangle_i^n)}{m(\triangle_i^n)}I_{\triangle_i^n}.$$

Then $A_n \geq CE(X_n^{\beta/2})$. The sequence $(X_{2^k}, k \geq 0)$ is a martingale with respect to the filtration $\mathcal{F}_{2^k}$. As a consequence (see, for instance, Theorem 3.3 in [2]), we have $\lim_{n\to\infty} X_{2^k} = X (m+\nu)$-a.e. Since $\nu \perp m$, $X = 0$ $m$-a.e. If $\lim_{k\to\infty} E(X_{2^k}^{\beta/2}) < \infty$, then $(X_{2^k}, k \geq 0)$ would be a uniformly integrable martingale and, hence, $X_{2^k} = E(X|\mathcal{F}_{2^k}) = 0$, which is a contradiction.

(ii) If $\alpha \in (0, \frac{1}{4})$, then

$$A_n = \sum_{i=1}^{n}\left(\int_0^{t_{i-1}}((t_i^n - s)^\alpha - (t_{i-1}^n - s)^\alpha)^2\, d\nu_s\right)^{\beta/2} + \sum_{i=1}^{n}\left(\int_{t_{i-1}}^{t_i^n}(t_i^n - s)^{2\alpha}\, d\nu_s\right)^{\beta/2}$$
$$= B_n + C_n.$$

For the term $C_n$ we have

$$C_n \leq \left(\frac{t}{n}\right)^{\alpha\beta}\sum_{i=1}^{n}(\nu(\triangle_i^n))^{\beta/2} = t^{\alpha\beta}\sum_{i=1}^{n}\frac{1}{n}(\nu(\triangle_i^n)n)^{\beta/2} = t^{\alpha\beta}E(X_n^{\beta/2}).$$

Since $E(X_n) = \nu([0,t]) < \infty$, $\frac{\beta}{2} < 1$, and $X_n \to 0$ a.e., we have $\lim_{n\to\infty} C_n = 0$. On the other hand,

$$B_n \leq \sum_{i=1}^{n}\left(\sum_{j=1}^{i-1}\int_{t_{j-1}}^{t_j}((t_i^n - s)^\alpha - (t_{i-1}^n - s)^\alpha)^2\, d\nu_s\right)^{\beta/2}$$
$$\leq \sum_{i=1}^{n}\left(\sum_{j=1}^{i-1}\left(\frac{t}{n}\right)^{2\alpha}(i^\alpha - (i-1)^\alpha)^2\nu(\triangle_j^n)\right)^{\beta/2}$$
$$\leq \sum_{i=1}^{n}\left(\sum_{j=1}^{i-1}\left(\frac{t}{n}\right)^{\alpha\beta}(i^\alpha - (i-1)^\alpha)^\beta\nu(\triangle_j^n)^{\beta/2}\right)$$
$$\leq \left(\frac{t}{n}\right)^{\alpha\beta}\sum_{i=1}^{n}(i^\alpha - (i-1)^\alpha)^\beta\sum_{j=1}^{n}\nu(\triangle_j^n)^{\beta/2}.$$

Notice that

$$\sum_{i=1}^{n}(i^{\alpha}-(i-1)^{\alpha})^{\beta} \leq C + \sum_{i=2}^{n}(i^{\alpha}-(i-1)^{\alpha})^{\beta}$$

$$\leq C + \sum_{i=2}^{n}(i-1)^{(\alpha-1)\beta} = C + O(n^{\alpha\beta-\beta+1}),$$

where $C > 0$. If $\alpha \in (0, \frac{1}{4})$, we have $\alpha\beta - \beta + 1 < 0$ and then $\sup_n \sum_{i=1}^{n}(i^{\alpha}-(i-1)^{\alpha})^{\beta} < \infty$. Then, similarly, $\lim_n A_n = 0$. □

**A.2. Transformations of Hölder continuous functions.** Let $\beta \in (0,1]$. We denote by $C^{\beta}([0,T])$ the set of Hölder continuous functions on $[0,T]$. For any function $f$ in $C^{\beta}([0,T])$ and any $0 \leq a < b \leq T$, we will write

(A.2) $$\|f\|_{\beta,a,b} = \sup_{a \leq s < t \leq b} \frac{|f(t) - f(s)|}{|t-s|^{\beta}}.$$

We also set $\|f\|_{\beta} = \|f\|_{\beta,0,T}$.

LEMMA A.4. *Suppose that $f \in C^{\beta}([0,T])$, and assume that $0 \leq a < b < v \leq T$. Let, $\gamma \geq 0$ and $\alpha + \beta \neq 0$. Then*

$$\left|\int_a^b s^{\gamma}(v-s)^{\alpha}\,df(s)\right| \leq \|f\|_{\beta}\left(2 + \left|\frac{\alpha}{\alpha+\beta}\right|\right)b^{\gamma}((v-b)^{\alpha+\beta} + (v-a)^{\alpha+\beta}).$$

PROOF. Suppose first $\gamma > 0$. Integrating by parts yields

$$\left|\int_a^b s^{\gamma}(v-s)^{\alpha}\,df(s)\right|$$

$$= \left|b^{\gamma}(v-b)^{\alpha}(f(b)-f(v)) - a^{\gamma}(v-a)^{\alpha}(f(a)-f(v))\right.$$

$$\left. - \int_a^b (f(s)-f(v))[s^{\gamma}(v-s)^{\alpha}]'\,ds\right|$$

$$\leq \|f\|_{\beta,a,v}\left(b^{\gamma}(v-b)^{\alpha+\beta} + a^{\gamma}(v-a)^{\alpha+\beta}\right.$$

$$\left. + \gamma\int_a^b (v-s)^{\alpha+\beta}s^{\gamma-1}\,ds + \alpha\int_a^b (v-s)^{\alpha+\beta-1}s^{\gamma}\,ds\right)$$

$$\leq \|f\|_{\beta,a,v}\left[b^{\gamma}(v-b)^{\alpha+\beta} + b^{\gamma}(v-a)^{\alpha+\beta}\right.$$

$$+ \max\{(v-a)^{\alpha+\beta}, (v-b)^{\alpha+\beta}\}(b^{\gamma}-a^{\gamma})$$

$$\left. + b^{\gamma}\left|\frac{\alpha}{\alpha+\beta}\right|((v-a)^{\alpha+\beta} - (v-b)^{\alpha+\beta})\right]$$

$$\leq \|f\|_{\beta,a,v}\left(2+\left|\frac{\alpha}{\alpha+\beta}\right|\right)b^{\gamma}((v-b)^{\alpha+\beta}+(v-a)^{\alpha+\beta}).$$

The case $\gamma = 0$ is proved in a similar way. $\square$

LEMMA A.5. *Suppose that* $f \in C^{\beta}([0,T])$, *and suppose* $\alpha < 0, \alpha + \beta > 0$. *Let* $g(t) = \int_0^t s^{\alpha} \, df(s)$. *Then,* $g \in C^{\alpha+\beta}([0,T])$, *and*

$$\|g\|_{\alpha+\beta} \leq \frac{\beta}{\alpha+\beta}\|f\|_{\beta}.$$

PROOF. Fix $0 \leq a < b \leq T$. Integrating by parts yields

$$|g(b) - g(a)| = \left|\int_a^b s^{\alpha} d[f(s) - f(a)]\right|$$

$$= \left|b^{\alpha}[f(b) - f(a)] + \alpha \int_a^b [f(s) - f(a)]s^{\alpha-1}\, ds\right|$$

$$\leq \|f\|_{\beta} b^{\alpha}|b-a|^{\beta} + |\alpha|\int_a^b |f(s) - f(a)|(s-a)^{\alpha-1}\, ds$$

$$\leq \|f\|_{\beta}\left(|b-a|^{\alpha+\beta} + |\alpha|\int_a^b (s-a)^{\alpha+\beta-1}\, ds\right)$$

$$\leq \|f\|_{\beta}\frac{\beta}{\alpha+\beta}|b-a|^{\alpha+\beta},$$

which give the desired result. $\square$

PROPOSITION A.6. *Fix* $\alpha \in (-\frac{1}{2}, \frac{1}{2})$ *and* $\beta \in (0,1]$ *such that* $0 < \alpha + \beta \leq 1$. *Suppose that* $f \in C^{\beta}([0,T])$, *and let* $g(t) = \int_0^t s^{\alpha}(t-s)^{\alpha}\, df_s$. *Then:*

1. *If* $\alpha > 0$, $g \in C^{\alpha+\beta}([0,T])$ *and for any* $0 \leq a < b \leq T$, *we have*

(A.3) $$|g(b) - g(a)| \leq C\|f\|_{\beta} b^{\alpha}(b-a)^{\alpha+\beta}.$$

2. *If* $\alpha < 0$ *and* $0 < 2\alpha + \beta \leq 1$, *then* $g \in C^{2\alpha+\beta}([0,T])$ *and*

$$|g(b) - g(a)| \leq C\|f\|_{\beta}(b-a)^{2\alpha+\beta}.$$

PROOF. We can write

$$g(b) - g(a) = \int_0^a s^{\alpha}((b-s)^{\alpha} - (a-s)^{\alpha})\, df_s + \int_a^b s^{\alpha}(b-s)^{\alpha}\, df_s$$

$$= \alpha \int_a^b \left(\int_0^a s^{\alpha}(v-s)^{\alpha-1}\, df_s\right) dv + \int_a^b s^{\alpha}(b-s)^{\alpha}\, df_s$$

$$= A + B.$$

If $\alpha > 0$, using Lemma A.4 yields

$$|A| \leq C\|f\|_\beta a^\alpha \int_a^b ((v-a)^{\alpha+\beta} + v^{\alpha+\beta})\, dv$$

$$= C\|f\|_\beta a^\alpha [(b-a)^{\alpha+\beta} + b^{\alpha+\beta} - a^{\alpha+\beta}]$$

and

$$|B| \leq C\|f\|_\beta b^\alpha (b-a)^{\alpha+\beta},$$

which implies (A.3) follows. On the other hand, if $\alpha < 0$, the function $h(t) = \int_0^t s^\alpha\, df_s$ is $(\alpha+\beta)$-Hölder continuous by Lemma A.5, and $\|h\|_{\alpha+\beta} \leq C\|f\|_\beta$. Then, applying Lemma A.4 to the function $h$, we obtain the estimates

$$|A| \leq \alpha \left| \int_a^b \left( \int_0^a (v-s)^{\alpha-1}\, dh_s \right) dv \right|$$

$$\leq C\|f\|_\beta \int_a^b [(v-a)^{2\alpha+\beta-1} + v^{2\alpha+\beta-1}]\, dv$$

$$\leq C\|f\|_\beta [(b-a)^{2\alpha+\beta} + b^{2\alpha+\beta} - a^{2\alpha+\beta}],$$

and

$$|B| \leq \left| \int_a^b (b-s)^\alpha\, dh_s \right| \leq C\|f\|_\beta (b-a)^{2\alpha+\beta}.$$

The proof is complete. $\square$

LEMMA A.7. *Fix $\alpha \in (-\frac{1}{2}, \frac{1}{2})$ and $\beta \in (0,1]$ such that $0 < \alpha + \beta \leq 1$ and $0 < 2\alpha + \beta \leq 1$. Suppose that $f \in C^\beta([0,T])$, and let $g(t) = \int_0^t s^\alpha (t-s)^\alpha\, df_s$. Set*

$$h(t) = \int_0^t u^{-\alpha-1} \left( \int_0^u (u-s)^{-\alpha}\, dg_s \right) du.$$

*Then for any $0 \leq a < b \leq T$, we have*

$$|h(b) - h(a)| \leq C\|f\|_\beta (b^\beta - a^\beta).$$

PROOF. We have

$$(A.4) \qquad |h(b) - h(a)| \leq \int_a^b u^{-\alpha-1} \left| \int_0^u (u-s)^{-\alpha}\, dg_s \right| du.$$

Suppose first that $\alpha < 0$. Then, $\|g\|_{2\alpha+\beta} \leq C\|f\|_\beta$, and Lemma A.4 yields

$$(A.5) \qquad \left| \int_0^u (u-s)^{-\alpha}\, dg_s \right| \leq C\|f\|_\beta u^{\alpha+\beta}.$$

Substituting (A.5) into (A.4) yields the results. In the case $\alpha > 0$, the Hölder norm $\|g\|_{\alpha+\beta}$ in an interval $[0, u]$ is bounded by $Cu^\alpha \|f\|_\beta$, and Lemma A.4 yields

$$\left| \int_0^u (u-s)^{-\alpha} \, dg_s \right| \leq C \|f\|_\beta u^{\beta+\alpha}.$$

This completes the proof of the lemma. □

**A.3. Existence of singular Hölder continuous distribution functions.** Let $0 < H < 1$ and $\rho > 1$. Suppose that $X = (X_t, t \geq 0)$ is a zero mean Gaussian process with stationary increments and a variance $\sigma^2(t) = E(X_t^2)$ given by

(A.6) $$\sigma^2(t) = \int_0^\infty (1 - \cos(xt)) g(x) \, dx,$$

where $g(x) = x^{-2H-1} \mathbf{1}_{[0,2)}(x) + (|\log x|^\rho x)^{-1} \mathbf{1}_{[2,\infty)}(x)$. If we replace $g(x)$ by $g_H(x) = x^{-2H-1}$ in equation (A.6), then the process $X$ is a fractional Brownian motion with Hurst parameter $H$. Taking into account that $g(x) \geq Cg_H(x)$ for some constant $C > 0$, it follows that the process $X$ satisfies the local nondeterminism property in some interval $(0, d)$ (see Theorem 4.1 in [1]).

The following lemma implies the existence of finite measures on the real line which are singular with respect to the Lebesgue measure, and whose distribution function is Hölder continuous of order $\gamma$, for any $\gamma < 1$ on any finite interval.

LEMMA A.8. *Let $X$ be the Gaussian process introduced above. Then, there exists a version of its local time $L(t, x)$, jointly continuous in $t$ and $x$, with the following properties:*

(i) *For each $x \in \mathbb{R}$ and $\gamma < 1$, $L(t, x)$ is Hölder continuous of order $\gamma$ with respect to $t$, on any finite interval.*

(ii) *$L(t, x)$ is a nondecreasing function of $t$.*

(iii) *For each $x \in \mathbb{R}$, the support of the measure $L(dt, x)$ is the set $\{s, X_s = x\}$, which has a Lebesgue measure $0$.*

PROOF. The function $\sigma^2$ satisfies

$$\sigma^2(t) \geq C |\log t^{-1}|^{-\alpha},$$

for some constant $C > 0$ and for $t \in (0, \frac{1}{2})$. Then, property (i) follows by Theorem 8.1 in [1]. From Theorem 6.4, page 11, in [3], it follows that for each $x \in \mathbb{R}$ the support of the measure $L(dt, x)$ is the set $\Lambda_x = \{s, X_s = x\}$. Finally, to show that $\Lambda_x$ has a Lebesgue measure $0$, we write

$$E \int_0^T \mathbf{1}_{\Lambda_x}(s) \, ds = \int_0^T E(\mathbf{1}_{X_s = x}) \, ds = 0,$$

which implies that $\int_0^T \mathbf{1}_{\Lambda_x}(s)\,ds = 0$ almost surely. This completes the proof of the lemma. $\square$

**Acknowledgments.** We would like to thank two anonymous referees for their helpful remarks.

## REFERENCES

[1] BERMAN, S. M. (1973/74). Local nondeterminism and local times of Gaussian processes. *Indiana Univ. Math. J.* **23** 69–94. MR0317397
[2] DURRETT, R. (2005). *Probability: Theory and Examples*, 3nd ed. Duxbury Press, Belmont, CA. MR1609153
[3] GEMAN, D. and HOROWITZ, J. (1980). Occupation densities. *Ann. Probab.* **8** 1–67. MR556414
[4] HU, Y. (2005). Integral transformations and anticipative calculus for fractional Brownian motions. *Mem. Amer. Math. Soc.* **175**. MR2130224
[5] KARATZAS, I. and SHREVE, S. E. (1991). *Brownian Motion and Stochastic Calculus*, 2nd ed. Springer, New York. MR1121940
[6] MANDELBROT, B. B. and VAN NESS, J. W. (1968). Fractional Brownian motions, fractional noises and applications. *SIAM Rev.* **10** 422–437. MR0242239
[7] MISHURA, J. and VALKEILA, E. (2007). An extension of the Lévy characterization to fractional Brownian motion. Preprint.
[8] NORROS, I., VALKEILA, E. and VIRTAMO, J. (1999). An elementary approach to a Girsanov formula and other analytical results on fractional Brownian motions. *Bernoulli* **5** 571–587. MR1704556
[9] NUALART, D. (2003). Stochastic integration with respect to fractional Brownian motion and applications. *Contemp. Math.* **336** 3–39. MR2037156
[10] ROGERS, L. C. G. (1997). Arbitrage with fractional Brownian motion. *Math. Finance* **7** 95–105. MR1434408
[11] SAMKO, S. G., KILBAS, A. A. and MARICHEV, O. I. (1993). *Fractional Integrals and Derivatives: Theory and Applications*. Gordon and Breach, Yverdon. MR1347689

DEPARTMENT OF MATHEMATICS
UNIVERSITY OF KANSAS
LAWRENCE, KANSAS 66045
USA
E-MAIL: hu@math.ku.edu